\documentclass[12pt]{article}
\usepackage{a4}
\usepackage{amssymb}
\usepackage{amsmath}
\usepackage{amsthm}
\usepackage{amstext}
\usepackage{amscd}
\usepackage{latexsym}
\usepackage{graphics}
\usepackage[OT2,T1]{fontenc}
\theoremstyle{plain}
\newtheorem{thm}{Theorem}[section]
\newtheorem{lemma}[thm]{Lemma}
\newtheorem{lem}[thm]{Lemma}
\newtheorem{cor}[thm]{Corollary}
\newtheorem{prop}[thm]{Proposition}
\theoremstyle{definition}
\newtheorem{rmk}[thm]{Remark}

\DeclareSymbolFont{cyrletters}{OT2}{wncyr}{m}{n}
\DeclareMathSymbol{\Sha}{\mathalpha}{cyrletters}{"58}

\newcommand{\F}{{\mathbb F}}

\renewcommand{\P}{{\mathbb P}}
\newcommand{\Q}{{\mathbb Q}}
\newcommand{\R}{{\mathbb R}}

\newcommand{\Z}{{\mathbb Z}}
\newcommand{\Gal}{{\rm Gal}}

\input{xy}
\xyoption{all}

\begin{document}

\title{The Chow group of zero-cycles on certain $\rm{Ch\hat{a}telet}$ surfaces over local fields}
\author{Supriya Pisolkar}
\date{}
\maketitle

\section{Introduction}\label{intro}
Let $\rm{K}$ be a finite extension of $\Q_p \ (p\ \ {\rm prime} )$. By a $\rm{Ch\hat{a}telet}$ surface $\rm{X}$ over $\rm{K}$ we mean a smooth projective surface $\rm{K}$-birational to a surface given by the equation:
\begin{equation} y^2-dz^2 = f(x)
\end{equation}
where $f(x)$ is a monic cubic polynomial in $x$ with coefficients in $\rm{K}$. Our main aim is to compute the Chow group $\rm{A_0(X)_0}$ of $0$-cycles of degree zero modulo rational equivalence on such surfaces. The case where $f(x)$ splits into three linear factors has been considered in \cite{dalawat1} and \cite{dalawat2}. In this paper we consider the remaining cases, in which $f(x)$ is either irreducible or of the form $x(x^2-e)$, where $e \in {\rm K^*}$ is not a square.\\

\noindent If $d \in {\rm K^*}^2$, then the $\rm{Ch\hat{a}telet}$  surface defined by the equation $y^2-dz^2=x(x^2-e)$ is ${\rm K}$-birational to $\P^2_{\rm K}$. In fact, in this case the function field of this surface is ${\rm K}(x,u)$, where $u = y+\sqrt{d}z$. By (\cite{coray}, Prop.~6.1), ${\rm A}_0(\P^2_{\rm K})_0 = 0$. Since ${\rm A_0(X)_0}$ is a birational invariant of a smooth projective geometrically integral surface (\cite{coray}, Prop.~6.3), we get that ${\rm A_0(X)_0}$ is zero. Thus, we may assume that $d \notin {\rm K}^{*2}$.\\

\noindent  The main results of this paper are as follows.

\begin{thm}\label{L=E}  Let $\rm{X}$ be a {\rm Ch\^atelet} surface given by the equation $y^2-dz^2=x(x^2-e)$. Let ${\rm L}={\rm K}(\sqrt d)$ and ${\rm E}={\rm K}(\sqrt e)$. If $\rm{L}$ and $\rm{E}$ are isomorphic extensions of ${\rm K}$, then ${\rm A_0(X)_0} = \{0\}$.

\end{thm}

\begin{thm}\label{podd}
Suppose that $p\neq 2$. Assume that the quadratic extensions ${\rm L}={\rm K}(\sqrt{d})$ and ${\rm E}={\rm K}(\sqrt{e})$ are not isomorphic. Then $\rm{A_0(X)_0}$ is isomorphic to $\Z/2\Z$.
\end{thm}

\begin{thm}\label{p2}
Suppose that $\rm{K}=\Q_2$. Assume that ${\rm L}={\rm K}(\sqrt{d})$ and ${\rm E}={\rm K}(\sqrt{e})$ are non-isomorphic quadratic extensions of ${\rm K}$. 
\begin{enumerate}
\item[{\rm (1)}] Supppose that $\rm{L/K}$ is unramified. Then the group $\rm{A_0(X)_0}$ is isomorphic to 
        \begin{enumerate}
                \item[{\rm (i)}] $\{0\}$ \hspace{6mm} if \hspace{3mm} $\upsilon_{\rm K}(e) \equiv 0 \ (\rm{mod}\ 4)$.
    \item[{\rm (ii)}]  $\Z/2\Z$ \hspace{3mm} if \hspace{3mm} $\upsilon_{\rm K}(e) \not\equiv 0 \ (\rm{mod}\ 4)$. 
  \end{enumerate}                          
\item[{\rm (2)}] Suppose that $\rm{L/K}$ is a ramified extension. Then ${\rm A_0(X)_0}$ is isomorphic to $\Z/2\Z$. 
\end{enumerate}
\end{thm}

\begin{thm}\label{thm:irreducible}
Let ${\rm X}$ be a Ch\^atelet surface {\rm K}-birational to $y^2-dz^2=f(x)$ where $f(x)$ is an irreducible monic cubic polynomial in $x$ with coefficients in ${\rm K}$. Then ${\rm A_0(X)_0} = \{0\}$
\end{thm}

\noindent {\bf Acknowledgement}: I am thankful to Prof.~C.~S.~Dalawat for his help. This problem was suggested to me by him. I am indebted to Jo\"{e}l Riou for painstakingly going through the first draft of this paper and for important suggestions. I am grateful to Amit Hogadi for his interest and stimulating discussions.  I thank Maneesh Thakur for his encouragement.


\section{The method of computation}\label{section:computation}

Let ${\rm X} = {\rm X}_{d,e}$ denote the Ch{\^a}telet surface corresponding to the equation $$y^2-dz^2=x(x^2-e)$$  Let $\pi$ be a uniformiser of $\rm{K}$. The change of variables $\ x = \pi^2x'$, $\ y = \pi^{3}y' $, $\ z =\pi^{3}z' $ gives us 
 $${\rm X}_{d,e'} :  {y'}^2-d{z'}^2 = x'({x'}^2-e')$$
 where $e'=\pi^{-4}e$. Thus it is enough to consider the cases $\ v_{\rm K}(e) = 0,1,2,3$. Moreover, using yet another transformation $z \mapsto \lambda z$ for a suitable $\lambda \in {\rm K}^*$, it is clear that we need only to consider the cases $v_{\rm K}(d)=0,1$.  \\ 
   
\noindent Let $\rm{CH}_0({\rm X})$ = Chow group of zero cycles on ${\rm X}$ modulo rational equivalence.
${\rm A_0(X)_0}$ = Ker $(\rm{CH}_0({\rm X})\stackrel{deg}{\longrightarrow}  \Z)$.\\

\noindent We now describe a method due to  Colliot-Th\'el\`ene and Sansuc \cite{Col} which reduces the calculation of ${\rm A_0(X)_0}$ to a purely number-theoretic question. \\

\noindent The surface ${\rm X}$ comes equipped with a morphism  $f: {\rm X} \rightarrow \P_1$ whose fibres are  conics. We denote by ${\rm O}$ the singular point of the fibre above $\infty$. By (\cite{coray} Th\'eor$\grave{\rm e}$me C), the map  $$\gamma : {\rm X}(\rm K)\rightarrow {\rm A_0(X)_0}, \ \ \gamma(Q) = Q-O, $$
\noindent is surjective. We also have a natural injection  (see \cite{Col}) $$\phi: {\rm A_0(X)_0} \rightarrow \rm{H}^1(\rm K, S(\rm{\overline{K}})$$  where ${\rm{\overline K}}$ is an algebraic closure of ${\rm K}$ and ${\rm S}$ is the ${\rm K}$-torus whose character group is the ${\rm Gal(\overline{K}/K)}$-module Pic($\overline{{\rm X}}$) where ${\rm {\overline X}} = {\rm X}\times_{\rm K}{\rm {\overline K}}$.
Thus with the following indentifications (see \cite{Col}) $$\rm{H}^1({\rm K}, S(\rm{\overline K})) \rightarrow \rm{K^*/N_{L/K}L^*} \times \rm{E^*/ N_{LE/E}LE^*} \rightarrow (\Z/2\Z)^2 $$ 
the calculation of ${\rm A_0(X)_0}$ reduces to computing the image of the composite map $${\rm X}(\rm K) \to {\rm A_0(X)_0} \to {\rm H^1(K, S({\overline K}))}\stackrel{\cong}\to (\Z/2\Z)^2.$$

\noindent As all the points in the same fibre of the map $f : {\rm X}(\rm K) \to \P^1(\rm K)$ are mutually equivalent $0$-cycles, what we have to compute is the image of the induced map $\chi: f({\rm X}(\rm K)) \to (\Z/2\Z)^2$. The subset $f({\rm X}(\rm K)) \subset \P^1(\rm K)$ is clearly equal to, 
 $${\rm M} = \{ x \in {\rm K^*} | x(x^2-e) \in \rm{N_{L/K}L^*} \} \cup \{ 0 \}. $$
\noindent The exact description of the  map $\chi: {\rm M} \to (\Z/2\Z)^2$ is given by (see \cite{coray}, \cite{san} )

 \[ \chi(x) = \left\{ \begin{array}{ll}(\overline{x}~,~ (x-\sqrt{e})^{-}) & \mbox{if $x \neq 0$ }\\
                                (-{\overline e}~,~(-\sqrt{e})^{-}) & \mbox{if $x=0$}, \end{array} \right. \] 

\noindent where the bar denotes the image in ${\rm K^*/N_{L/K}L^*}$ and ${\rm E^*/N_{LE/E}LE^*}$ respectively, both these quotients being identified with $\Z/2\Z$. By using this map $\chi$ we will now prove Theorem \ref{L=E}.
\begin{proof}[Proof of Theorem \ref{L=E}] By the above method, to show that ${\rm A_0(X)_0} = \{0\}$, it is enough to show that $\chi({\rm M}) = \{0\}$. As ${\rm L}$ and ${\rm E}$ are isomorphic, the extension $\rm {LE/E}$ is trivial. Thus the group ${\rm E^*}/{\rm N_{LE/E}LE^*}$ is trivial. Therefore for any $x \in {\rm M}$ we get $\chi(x) = (\overline{x},0)$.  Since ${\rm N_{L/K}L^*} = {\rm N_{E/K}E^*}$, $-e \in {\rm N_{L/K}L^*}$. Thus $\chi(0)=(0,0)$.  Now let $x\in {\rm M}\backslash \{0\}$. Note that $x^2-e  \in {\rm N_{E/K}E^*}={\rm N_{L/K}L^*}$. This, together with the fact that $x(x^2-e)\in {\rm N_{L/K}L^*}$, implies that  $x \in {\rm N_{L/K}L^*}$ and $\chi(x) = (0,0)$.

\end{proof}

\noindent Before  proving Theorem \ref{podd} we observe that $\chi(\rm M)$ is contained in the diagonal subgroup of $(\Z/2\Z)^2$ when ${\rm L}$ and ${\rm E}$ are non-isomorphic.

\begin{lemma} \label{lem2}
Let $\ \rm{L/K} = {\rm K}(\sqrt{d})$ and ${\rm E/K} = {\rm K}(\sqrt{e})$ be non-isomorphic quadratic extensions. Then $\chi(M)$ is  contained in the diagonal subgroup of $\ \Z/2\Z \times \Z/2\Z$. In particular, ${\rm A_0(X)_0}$ is either $\{0\}$ or $\Z/2\Z$.
\end{lemma}

\begin{proof}  By class field theory (see \cite{serre2}, p.~212) we  have the commutative diagram\\
$$
\xymatrix{
\rm{E^*/N_{LE/E}LE^*}\ar[r]^{\rm{N_{E/K}}}\ar[d]_{\rm{rec}}   &  \rm{K^*/N_{L/K}L^*}\ar[d]^{\rm{rec}}   \\
\rm{Gal(LE/E)}\ar[r]^{\cong}    &          \rm{Gal(L/K)}
}
$$
where the vertical maps are isomorphisms. The map ${\rm Gal(LE/E)} \to {\rm Gal(L/K)}$ is  an ismorphism  since $\rm L$ and $\rm E$ are linearly disjoint. Thus $\rm{E^*/N_{LE/E}LE^*} \to {\rm K^*/N_{L/K}L^*}$ is an isomorphism, i.e. an element $t\in {\rm E}^*$ belongs to $\ {\rm N_{LE/E}LE^*} \ $ if and only if $\ {\rm N_{E/K}}(t)\ $ belongs to ${\rm N_{L/K}L^*}$. Therefore, for any $ x \in \rm{M}\backslash\{0\}$, $ x \in {\rm N_{L/K}L^*}$ if and only if $x-\sqrt{e} \in {\rm N_{LE/E}LE^*}$. Thus $\chi(x) = (0,0)$ or $(1,1)$. Similarly, $\chi(0) = (-e, -\sqrt{e})=(0,0)$ or $(1,1)$  depending upon whether $-e \in \rm{N_{L/K}L^*}$ or $-e \not\in \rm{N_{L/K}L^*}$. This shows that  $\chi(\rm{M})$ is contained in the diagonal subgroup of $\Z/2\Z \times \Z/2\Z$. 

\end{proof}

\begin{cor}\label{cor:diagonal}  The group $\rm{A_0(X)_0}\cong \Z/2\Z ~$ if and only if at least one of the following condition holds.
\begin{enumerate}
\item[(i)]$~-e \not\in \rm{N_{L/K}L^*}$.
\item[(ii)] There exists $x \in {\rm K^*}$, such that $x \notin {\rm N_{L/K}L^*}$ and $x^2-e \notin {\rm N_{L/K}L^*}$.
\end{enumerate}
\end{cor}

\section{Proof of Theorem \ref{podd}}
\noindent In this section we prove Theorem \ref{podd}. Throughout this section let $p$ denote an odd prime and let $\rm{K}$ denote a finite extension of $\Q_p$. 

\begin{lemma}\label{lem3}
Let $\rm{F/K}$ be a quadratic extension of $\rm{K}$. 
\begin{enumerate}
\item[(i)] If $~\rm{F/K}$ is unramified, then an element $x \in \rm{K}^*$ belongs to $\rm{N_{F/K}F^*}$ if and only if $v_{\rm K}(x)$ is even.
\item[(ii)] If $~{\rm F/K}$ is ramified, $\pi_{\rm F}$ is a uniformiser of ${\rm F}$ and $\pi_{\rm K}={\rm N_{F/K}}(\pi_{\rm F})$, then $x\in \rm{K^*}$ belongs to $\rm{N_{F/K}F^*}$ if and only if $x/\pi_{\rm K}^{v_{\rm K}(x)}$ is a square.
\end{enumerate}
\end{lemma}
\begin{proof}
\noindent (i) It is easy to see that ${\rm N_{F/K}F^*}$ is contained in the subgroup of elements of even valuation. Since both these subgroups are of index two in ${\rm K^*}$, they are equal. \\
\noindent (ii) Let ${\rm N'}$ be the subgroup of all elements $x\in {\rm K}^*$ such that $x/\pi_{\rm K}^{v_{\rm K}(x)}$ is a square. Using the fact that $\pi_{\rm K}$ belongs to ${\rm N_{F/K}F^*}$, it is clear that we have ${\rm N'} \subset {\rm N_{F/K}F^*}$. Since  ${\rm N'}$ and ${\rm N_{F/K}F^*}$ are index-two subgroups of $\rm K^*$, they must be equal. 

\end{proof}

\begin{lemma}\label{lem4} Let ${\rm K}$ be a finite extension of $\Q_p$ where $p$ is an odd prime.  Suppose that $e \in {\rm K^*}$ is not a square. Then $\ \rm{E} = \rm{K}(\sqrt{\it e})\ $ is a ramified extension of $\ \rm{K}\ $ if and only if $\ v_{\rm K}(e) $ is odd. 
\end{lemma}
\begin{proof} If $v_{\rm K}(e)$ is odd, then we make the reduction to the case where $v_{\rm K}(e)=1$ by multplying $e$ by a square. It is clear that ${\rm E}$ is ramified when $v_{\rm K}(e)= 1$. Now suppose $v_{\rm K}(e)$ is even.  We may assume that $e$ is a unit by modifying $e$ by a square. Since $p$ is odd, the polynomial $T^2-e$ is separable over the residue field, and hence irreducible in the residue field. Thus ${\rm E/K}$ is unramified in this case.   

\end{proof}

\begin{proof}[Proof of Theorem \ref{podd}]
We split the proof into following two cases. To show that $\rm{A_0(X)_0}$ is isomorphic to $\Z/2\Z$, by lemma \ref{lem2}, it suffices to show that there exists an $x \in {\rm M}$ such that $\chi(x)=(1,1)$ in each of these cases.\\

\noindent Case~($1$) : $\rm{E}$ is unramified over $\rm{K}$. \\
\noindent In this case ${\rm L/K}$, being non-isomorphic to ${\rm E/K}$, is a ramified extension. If $-e \notin \rm{N_{L/K}L^*}$ then $\chi(0) = (1,1)$ and we are done. Suppose that $-e \in \rm{N_{L/K}L^*}$. Let us first show that $e \notin {\rm N_{L/K}L^*}$. Indeed, as ${\rm E} = {\rm K}(\sqrt{e})$ is the unramified quadratic extension, $v_{\rm K}(e)$ is even by Lemma \ref{lem4} and the unit $e\pi^{-v_{\rm K}(e)}$ is not a square for any uniformiser $\pi$ of ${\rm K}$. Lemma \ref{lem3}(ii) now implies that $e \notin \rm{N_{L/K}L^*}$. Hence $-1 \not\in \rm{N_{L/K}L^*}$; in particular, $-1$ is not a square.\\
\noindent Write $\ -1 = x^2-ey^2\ $ for some $x, y \in {\rm K}$. This is possible because $-1 \in {\rm N_{E/K}E^*}$ by Lemma \ref{lem3}(i). Clearly $\ x \neq 0 $ and $ y \neq 0 $, because neither $-1$ nor $e$ is a square.\\
\noindent Put $\ \alpha = x/y $. Replacing $\alpha$ by $-\alpha$ if necessary, we may assume that $\ \alpha \not\in \rm{N_{L/K}L^*}$. Moreover, $\alpha^2-e \notin {\rm N_{L/K}L^*}$ because $\alpha^2-e = -1/y^2$  and $-1 \notin {\rm N_{L/K}L^*}$.  It follows that $\alpha(\alpha^2-e) \in {\rm N_{L/K}L^*}$. Thus $\ \alpha \in {\rm M}$ and as $\chi(\alpha) =(1,1)$ we are done.\\


\noindent Case~($2$) : $\rm{E}$ is ramified over ${\rm K}$. \\
\noindent We will show that $\chi(0) = (1,1)$. As $\ \rm{E/K}\ $ and $\rm{L/K}$ are quadratic extensions, by local class field theory, their norm subgroups $\rm{N_{E/K}E^*}$ and $\rm{N_{L/K}L^*}$ are two index-two subgroups of $\rm{K^*}$. These two subgroups are not equal as $\rm{E}$ and $\rm{L}$ are not isomorphic. Then, their intersection must be an index-four subgroup of $\rm{K^*}$. Since $(\rm{K^*})^2 \subset \rm{N_{E/K}E^*} \cap \rm{N_{L/K}L^*}$ is also an index four subgroup of $\rm{K^*}$, we get $(\rm{K^*})^2 = {\rm N_{E/K}E^* \cap N_{L/K}L^*}$. We know that $\ -e  = {\rm N_{E/K}}(\sqrt{e}) \in {\rm N_{E/K}E^*} $. If $-e \in {\rm N_{L/K}L^*}$ then $-e$ would be a square. This contradicts the fact (by Lemma \ref{lem4}) that  $v_{\rm K}(e) = v_{\rm K}(-e)$ is odd. Thus $-e \not\in {\rm N_{L/K}L^*}$ and $\chi(0) = (1,1)$.

\end{proof}

\section{Preliminaries on Hilbert symbols}\label{hs}
\noindent In this section we review the notion of Hilbert symbol as given in \cite{serre1}.\\

\noindent Let $\rm{K}$ be a field. Let $a,b\in {\rm K^*}$. We put $(a,b)_{\rm K}= 1$ if $ax^2+by^2=1$ has a solution 
$(x,y)\in \rm{K^2}$ and $(a,b)_{\rm K}= -1$ otherwise. 
The number $(a,b)_{\rm K}$ is called the {\bf Hilbert symbol} of a and b relative to the field $\rm{K}$. It can be shown that the Hilbert symbol has the following properties.  
\begin{enumerate}
\item $(a,b)_{\rm K}=(b,a)_{\rm K}$ (Symmetry).
\item $(ab,c)_{\rm K}= (a,c)_{\rm K}\cdot (b,c)_{\rm K}$ (Bilinearity). 
\item $(a^2,b)_{\rm K}=1$.
\item $(a,1-a)_{\rm K}=1$ for $a \neq 1$.
\end{enumerate}

\noindent Thus the Hilbert symbol can be thought of as a symmetric bilinear form on the $\F_2$-vector space ${\rm K}^*/{\rm K}^{*2}$ with values in the group $\{1,-1\}$. \\

\noindent The following proposition gives an equivalent definition of the Hilbert symbol.  

\begin{prop} \label{normhilbert} {\rm (\cite{serre1}, p.~19)}
 Let $a,b \in \rm{K^*}$ and let ${\rm K}_b = {\rm K}(\sqrt{b})$. 
Then $(a,b)_{\rm K}= 1$ if and only if $a\in {\rm N}_{{\rm K}_b/{\rm K}}({\rm K}_{b}^*)$.
\end{prop}

\noindent When ${\rm K}=\Q_p$ we denote the Hilbert symbol by $(a,b)_p$ instead of $(a,b)_{\Q_p}$.\\

\noindent We now describe a formula for calculating the Hilbert symbol when ${\rm K}=\Q_2$. Let $\Z_2^*$ be the group of units of $\Z_2$ and ${\rm U}_3=1+8\Z_2$. Let $\varepsilon, \omega  :\Z_2^* \longrightarrow \Z/2\Z$ be the homomorphisms given by

$$\varepsilon(z) = \frac{z-1}{2}(\text{mod}\ 2 ) , \ \ \ \ \omega(z) = \frac{z^2-1}{8} (\text{mod}\ 2).$$ 
  
\noindent Put $a =  2^{\alpha}u$,$ b=2^{\beta}v$ where $u$ and $v$ are units. Then, according to {\rm (\cite{serre1}, p.~20)},  
$$(a,b)_2 = (-1)^{\varepsilon(u)\varepsilon(v)+ \omega(v)\alpha+ \omega(u)\beta}.$$

\noindent Using the above formula for the $2$-adic Hilbert symbol, we prove the following lemma which will be used in the proof of Theorem \ref{p2}.

\begin{lem}\label{e1d1} Let ${\rm K=\Q_2}$. Let $e,d,{\rm L,E}$ be as in Theorem \ref{p2}. Assume that $v_{\rm K}(d)=1$. Then, 
\begin{enumerate}
\item[(i)]If $v_{\rm K}(e)=1$, then at least one of the elements $-1,\ 1-e, \  e$ does not belong to ${\rm N_{L/K}L^*}$ 
\item[(ii)] If $v_{\rm K}(e)=3$, then at least one of the elements $-1,\ (1-e/4),\  e$ does not belong to ${\rm N_{L/K}L^*}$. 
\end{enumerate}
\end{lem}
\begin{proof}
\noindent $(i)$ Assume that all three elements $-1,1-e,e$ belong to ${\rm N_{L/K}L^*}$. Then by Prop. \ref{normhilbert}, $(-1,d)_2=(e,d)_2=(1-e,d)_2=1$. Let $e=2u$ and $d=2v$ where $u, v \in \Z_2^*$. We have, $$(-1,d)_2=(-1,2v)_2=(-1,2)_2(-1,v)_2= (-1)^{\varepsilon (v)}.$$ 
As $(-1,d)_2=1$ by assumption, we have \\
\begin{equation} \label{eqn:v14}
\varepsilon(v)= 0 \ \text{and thus} \ v \equiv 1 \ ({\rm mod}\ 4). 
\end{equation} 

\noindent Also, $$(1-e,d)_2=(1-2u,2v)_2=(-1)^{\varepsilon(1-2u)\varepsilon(v)+ \omega(1-2u)}$$ \\
As $\varepsilon(v) = 0 $ and by assumption  $(1-e,d)_2=1$, we get the following\\
\begin{equation}\label{eqn:u14}
 \omega(1-2u) = 0 \ \text{and} \ u \equiv 1 ({\rm mod}\ 4).
\end{equation}
\noindent Now,
$$(e,d)_2=(2u,2v)_2= (-1)^{\varepsilon(u)\varepsilon(v)+ \omega(v)+ \omega(u)}$$
As  $\varepsilon(v)= 0$ and $(e,d)_2=1$  by assumption, we get $\omega(u)=\omega(v)$. Since both $u,v$ are congruent to $1$ modulo $4$ by (\ref{eqn:v14}) and (\ref{eqn:u14}),  one can check that $u\equiv v ({\rm mod}\ 8)$. As any element of  ${\rm U}_3$ is a square, we get that $u$ and $v$ differ by a square unit. Thus the extensions ${\rm L}=\Q_2(\sqrt{2v})$ and ${\rm E}=\Q_2(\sqrt{2u})$ are isomorphic. This contradicts the hypothesis that ${\rm L}$ and ${\rm E}$ are non-isomorphic extensions of ${\rm K}$.  \\

\noindent (ii) Since the proof of this case is similar to the one above, we only give a sketch. Assume that all three elements $-1,(1-e/4),e$ belong to ${\rm N_{L/K}L^*}$. Let $e=2^3u$ and $d=2v$ where $u,v\in \Z_2^*$. As in $(i)$, using $(-1,d)_2=1$ we get $\varepsilon(v)=0$ and thus $v \equiv 1 \ {\rm mod}\ 4 $. Similarly $(1-e/4,d)_2=1$ gives $\omega(1-2u)=1$ and thus $u \equiv 1 \ {\rm mod}\ 4 $. By properties $2.$ and $3.$ of the Hilbert symbol mentioned earlier, we deduce that 
$$ (e,d)_2 = (2^3u,2v)_2 = (2u,2v)_2$$ 
Thus $(e,d)_2=1$ gives $\omega(u)=\omega(v)=1$. As in $(i)$, we arrive at a contradiction by showing that ${\rm L}$ and ${\rm E}$ are isomorphic extension of $\Q_2$.  

\end{proof}
\section{Some results on ramified quadratic extension}\label{s}

Throughout this section $\rm{K}$ will denote a finite extension of $\Q_2$. Let $\rm{L/K}$ be a ramified quadratic extension. Let $k$ be the residue field of $\rm{K}$. Since ${\rm L/K}$ is totally ramified, $k$ is also the residue field of ${\rm L}$. For $i\geq 0$, let ${\rm U}_{i,{\rm L}} = \{ x \in {\rm U_L} | v_{\rm L}(1-x) \geq i\}$. The subgroups $\{ {\rm U}_{i,{\rm L}}\}_{i \geq 0}$ define a decreasing filtration of $\rm{U_L}$. Similarly we define ${\rm U}_{i,{\rm K}}$.

\begin{thm}{ \rm(\cite{fesenko}, III.1.4) } \label{fesenko1}
Let $\pi_{\rm L}$ be a uniformiser of ${\rm L}$. Let $\sigma$ be the generator of $\Gal(\rm{L/K})$. 
Then $\frac{\sigma(\pi_{\rm L})}{\pi_{\rm L}} \in  {\rm U_{1,L}}$. 
Further, if $s$ is the largest integer such that 
$\frac{\sigma(\pi_{\rm L})}{\pi_{\rm L}}
\in {\rm U}_{s,\rm{L}}$,  
then $s$ is independent of the uniformiser $\pi_{\rm L}$. 
\end{thm}

\noindent Thus the integer $s$ defined above depends only on the extension $\rm{L/K}$. Therefore we will denote it by $s({\rm L/K})$.

\begin{thm}{ \rm (\cite{fesenko},~III.2.3) }\label{fesenko2}
Let $e_{\rm K}$ be the ramification index of ${\rm K}$ over $\Q_2$. Then $s({\rm L/K}) \leq 2e_{\rm K}$. 
\end{thm}

\begin{rmk}\label{rmk:ex}
When ${\rm K}=\Q_2$, $\rm K$ has six ramified quadratic extensions, namely ${\rm K(\sqrt{-1}), K(\sqrt{-5}), K(\sqrt{\pm 2}), K(\sqrt{\pm 10})}$. For the first two extensions $s=1$ and for the remaining extensions $s=2$.\end{rmk}
\begin{thm}{ \rm (\cite{fesenko},III.1.5)}\label{fesenko3}
Let ${\rm L,K}$ be as above. Assume that $k=\F_2$. Let $s = s({\rm L/K})$. Choose a uniformiser $\pi_{\rm L}$ of ${\rm L}$. Let $\pi_{\rm K}= {\rm N_{L/K}(\pi_L)}$. Define $$\lambda_{i,{\rm L}}:{\rm U}_{i,{\rm L}}\to \F_2 \hspace{2mm};\hspace{3mm} \lambda_{i,{\rm L}}(1+\theta \pi_{\rm L}^i)= \overline{\theta}$$
Similarly define $\lambda_{i,{\rm K}}$ using the uniformiser $\pi_{\rm K}$. Then 
\begin{enumerate}
\item The following diagrams commute.
$$
\xymatrix{
{\rm U}_{i,{\rm L}}\ar[r]^{\lambda_{i,{\rm L}}}\ar[d]_{\rm{N_{L/K}}}            &    \F_2  
    \ar[d]^{\rm{Id}} & {\rm if}\  1 \leq i < s \\
{\rm U}_{i,{\rm K}}\ar[r]^{\lambda_{i,{\rm K}}}  &  \F_2  & \\
{\rm U}_{s,{\rm L}}\ar[r]^{\lambda_{s,{\rm L}}}\ar[d]_{\rm{N_{L/K}}}    & \F_2 \ar[d]^0     &               \\
{\rm U}_{s,{\rm K}}\ar[r]^{\lambda_{s,{\rm K}}} & \F_2  & \\
{\rm U}_{s+2i,{\rm L}}\ar[r]^{\lambda_{s+2i,{\rm L}}}\ar[d]_{\rm{N_{L/K}}}              &               \F_2       \ar[d]^{\rm Id} & {\rm if}\ i>0 \\
{\rm U}_{s+i,{\rm K}}\ar[r]^{\lambda_{s+i,{\rm K}}}  &   \F_2  & 
}
$$
\item ${\rm N_{L/K}}({\rm U}_{s+i,{\rm L}})={\rm N_{L/K}}({\rm U}_{s+i+1,{\rm L}})$ for $i>0$, $p\nmid i$.
\item ${\rm N_{L/K}}({\rm U}_{s+1,{\rm L}})={\rm U}_{s+1,{\rm K}}$.
\end{enumerate}
\end{thm}

\noindent The following corollary is an easy consequence of the above theorem. However part $(iii)$ of the corollary will play an important role in the proof of Theorem \ref{p2}. 

\begin{cor}\label{important}
 With the notation as above, ${\rm N_{L/K}}({\rm U}_{i,{\rm L}})\subset {\rm U}_{i,{\rm K}}$ for $1 \leq i \leq s+1$ . Thus for $i \leq s$, we have induced maps
$${\rm N}^i_{\rm L/K} : \frac{ {\rm U}_{i,{\rm L}}}{{\rm U}_{i+1,{\rm L}}} \to \frac{ {\rm U}_{i,{\rm K}}}{{\rm U}_{i+1,{\rm K}}}.$$ Further,
\begin{enumerate}
\item[$(i)$] ${\rm N}^i_{\rm L/K}$ is an isomorphism for $1 \leq i < s$.
\item[$(ii)$] ${\rm N}^s_{\rm L/K}$ is the zero map. 
\item[$(iii)$] $\gamma \in {\rm U}_{s,{\rm K}}\backslash {\rm U}_{s+1,{\rm K}}\implies \gamma \notin {\rm N_{L/K}L}^*$.
\end{enumerate}
\end{cor}
\begin{proof}
$(1)$ and $(2)$ are immediate consequences of the first two commutative diagrams in Theorem \ref{fesenko3} and the fact that ${\rm U}_{i+1,{\rm L}}$ is the kernel of $\lambda_{i,{\rm L}}$. Suppose that $\gamma \in {\rm U}_{s,{\rm K}}\backslash {\rm U}_{s+1,{\rm K}}$. This gives that $\lambda_{s,{\rm K}}(\gamma)\neq 0$ and thus the second commutative diagram in \ref{fesenko3} tells us that $\gamma \notin {\rm N_{L/K}}({\rm U}_{s,{\rm L}})$. Now we want to show that for $i \neq s$ and for any $x\in {\rm U}_{i,{\rm L}}\backslash {\rm U}_{i+1,{\rm L}}$, ${\rm N}_{\rm L/K}(x) \neq \gamma$. For $i<s$ this is easy to see by (i) above. For $i > s$,
 \ref{fesenko3}(3) implies that ${\rm N}_{\rm L/K}(x) \in {\rm U}_{s+1,{\rm K}}$. Thus $\gamma \neq  {\rm N_{L/K}}(x)$.

\end{proof}


\section{The proof of Theorem \ref{p2}}\label{q2}

\noindent We first state the following lemma from \cite{serre2}.

\begin{lem} \label{lemma:square} {\rm (\cite{serre2}, p.~212)}
Let ${\rm K}/\Q_2$ be a finite extension.
Let $e_{\rm K}$ be the ramification index of ${\rm K}/\Q_2$.
For every $n > e_{\rm K}$ the map ${\rm U}_{n,{\rm K}} \stackrel{ x\mapsto x^2}{\longrightarrow} {\rm U}_{n+e_{\rm K},{\rm K}}$ is an isomorphism of $\Z_2$-modules.
\end{lem}

\noindent Now the proof of Theorem \ref{p2} will occupy the rest of this section. Henceforth ${\rm K}=\Q_2$. For simplicity of notation we write $\Z_2^*$ instead of ${\rm U}_{\Q_2}$, ${\rm U}_n$ instead of ${\rm U}_{n,\Q_2}$ and $v$ instead of $v_{\Q_2}$. For any field extension ${\rm F}/\Q_2$, we will write ${\rm N}({\rm F}^*)$ instead of ${\rm N}_{{\rm F}/\Q_2}({\rm F}^*)$.

\begin{proof}[Proof of Theorem \ref{p2}]Recall from section \ref{section:computation} that it is enough to consider the cases $v(d) = 0,1$ and $v(e)=0,1,2$ or $3$.  \\
\noindent {\bf (1)  $\rm{L}$ is unramified over $\Q_2$. } 

\noindent {\bf Case}: $v(e) = 0$ : To show that $\ \rm{A_0(X)_0}\ $ is zero, we have to show that $\chi(x) = 0$ for every $x \in {\rm M}$. 
If $x \in {\rm M}\backslash \{0\}$ such that $v(x) > 0$ then $$v(x^2-e) = {\rm min}\{2v(x),v(e)\} = 0.$$
If $x \in {\rm M}\backslash \{0\}$ such that $v(x) < 0$ then $$v(x^2-e) = {\rm min}\{2v(x),v(e)\} = 2v(x).$$
\noindent Thus whenever $v(x)\neq 0$, $v(x^2-e)$ is even and hence by Lemma \ref{lem3}, $x^2-e \in {\rm N(L^*)}$. Thus for any $x \in {\rm M}\backslash\{0\}$ with $v(x)\neq 0$, $\chi(x) = (0,0)$. If $x \in \rm{M}$ such that $v(x) = 0$ then $ x \in \rm{N(L^*)}$ and thus $\chi(x) = (0,0)$. We now claim that $\chi(0)=(0,0)$. For this we need to show that $-e \in \rm{N(L^*)}$. But this is clear since $v(-e) = 0$ and ${\rm L}$ is the unramified extension of $\Q_2$. This proves \ref{p2}(1)(i). \\

\noindent Now we prove \ref{p2}(1)(ii).\\

\noindent {\bf Case}: $v(e)= 1$ or $3$ : Choose $ x \in \Q_2^* $ such that $v(x)$ is odd and $2v(x) > v(e)$. Then by Lemma \ref{lem3}, $x \not\in {\rm N(L^*)}\ $ and $\ x^2-e \not\in {\rm N(L^*)}$. Therefore $x \in \rm{M}$ and $\chi(x) = (1,1)$ i.e ${\rm A_0(X)_0}$ is isomorphic to $\Z/2\Z$.\\ 

\noindent {\bf Case}: $v(e) = 2$ :  We claim that $\chi(2)=(1,1)$. Let $\beta = e/4$. Since $e$ is not a square in $\Q_2$, $\beta$ cannot be a square and thus $\beta \notin \rm{U}_3$. We now show that $\beta \notin {\rm U_2}$. Note that,  an element $\alpha \in {\rm U_2 \backslash U_3}$, can be written as $\alpha = 5.\gamma^2$ for some $\gamma \in {\rm U_1}$. Thus, $\Q_2({\sqrt 5}) = \Q_2({\sqrt \alpha})$, where $\Q_2(\sqrt{5})$ is the uramified quadratic extension of $\Q_2$. This implies that $\beta \notin {\rm U}_2$, since $\Q_2(\sqrt{\beta})= {\rm E}$ is a ramified quadratic extension of $\Q_2$.  Thus $\beta \in {\rm U_1\backslash U_2}$ and  $v(2^2-e)= v(2^2(1-\beta))= 3$. By Lemma \ref{lem3}(i),  $2^2-e \not\in {\rm N(L^*)}$ and similarly $2 \notin {\rm N(L^*)}$. Therefore $2 \in {\rm M} $ and $\chi(2) = (1,1)$ i.e  ${\rm A_0(X)_0}$ is isomorphic to $\Z/2\Z$. 

\noindent The proof of this case completes the proof of the Theorem \ref{p2}(1).

\vspace{3mm}
\noindent {\bf $(\rm{2})$ $\rm{L}$ is ramified over $\Q_2$} \\
\noindent We first prove the following lemma which will be crucially used in the proof.

\begin{lem}\label{-1e} Suppose that ${\rm L/\Q_2}$ is a ramified extension and either $-1$ or $e$ does not belong to ${\rm N_{L/K}L^*}$. Then $\rm{A_0(X)_0} \cong \Z/2\Z$.
\end{lem}
\begin{proof}
Suppose that $-1 \not\in {\rm N(L^*)}$. Since $\rm{L}$ and $\rm{E}$ are non-isomorphic quadratic extensions, by class field theory we can choose $\alpha \in \rm{N(E^*)}$ such that $\alpha \notin \rm{N(L^*)}$. Thus $\alpha = x^2-ey^2$ for some $x,y \in \Q_2$. Since $\alpha \notin \rm{N(L^*)}$, $y \neq 0$. If $x = 0$ then $-e \notin \rm{N(L^*)} \implies \chi(0)=(1,1) \implies \rm{A_0(X)_0} \cong \Z/2\Z$. Thus we assume that $x \neq 0$. 
Now $$\left(\frac{x}{y}\right)^2 - e = \frac{\alpha}{y^2} \notin \rm{N(L^*)}$$ If $\frac{x}{y} \notin \rm{N(L^*)}$ then $\frac{x}{y} \in \rm{M}$ and $ \chi(\frac{x}{y}) = (1,1)$. Otherwise since $-1 \notin \rm{N(L^*)}$, $-\frac{x}{y} \in \rm{M}$ and $\chi(-\frac{x}{y}) = (1,1)$. \\
\noindent Now suppose $e \not\in {\rm N(L^*)}$. We may assume $-1 \in \rm{N(L^*)}$ since otherwise we are done by the above case. Thus $-e \notin \rm{N(L^*)}$. Hence $\chi(0)=(1,1)$. 

\end{proof}

\noindent Now we prove the Theorem \ref{p2}(2). The proof splits into two parts depending upon the invariant $s({\rm L}/\Q_2)$ of the field extension. Since ${\rm L}/\Q_2$ is a ramified quadratic extension, $s({\rm L}/\Q_2)= 1$ or $2$, by Theorem \ref{fesenko2}. We want to show that ${\rm A_0(X)_0} \cong  \Z/2\Z$, i.e., to show that there exists an element $x \in \rm M$ such that $\chi(x) = (1,1)$ in each of these cases.\\  

\noindent \underline{{\bf Case $s({\rm L/K})=1$}}:  It is clear  that $-1 \in {\rm U}_{1} \backslash {\rm U}_{2}$. Thus by Cor.~\ref{important}(iii), $-1 \not\in \rm{N(L^*)}$. Thus ${\rm A_0(X)_0} \cong \Z/2\Z$ by Lemma \ref{-1e}.


\vspace{3mm}
\noindent \underline{{\bf Case $s({\rm L/K})=2$}}:  We may assume that $-1$ and $e$ belong to $\rm{N(L^*)}$, since otherwise  we are done by  Lemma \ref{-1e}. 

\noindent {\bf Case $v(e)=0$}:\\
Choose a uniformiser $\pi$ of $\Q_2$ which does not belong to  $\rm{N(L^*)}$.  Put $x = \pi e$. Since we have assumed that $e \in \rm{N(L^*)}$, we get $x \not\in \rm{N(L^*)}$. Write $x^2-e = -e(1-\pi^2e)$. Then $1-\pi^2e \in {\rm U}_2 \backslash  {\rm U}_3$ and by Cor.~\ref{important}(iii) we get $1-\pi^2e \notin {\rm N(L^*)}$. Since $-e \in {\rm N(L^*)}$, $x^2-e \not\in \rm{N(L^*)}$. Therefore $x \in {\rm M}$ and $\chi(x)=(1,1)$. \\

\noindent {\bf Case $v(e)=1$}: 
Since $s({\rm L}/\Q_2) = 2$, by remark \ref{rmk:ex} we know that $v(d) = 1$. In this case, by Lemma \ref{e1d1}(i), at least one of the elements $-1, e , 1-e$ does not belong to ${\rm N(L^*)}$. By our assumption, $-1,e \in {\rm N(L^*)}$. Thus $-e \in {\rm N(L^*)}$ and by Lemma \ref{e1d1}(i), $1-e \notin {\rm N(L^*)}$.\\
Put $x=eu$ where $u$ is any unit which is not in $\rm{N(L^*)}$. Thus $x \notin {\rm N(L^*)}$. Then $x^2-e = -e(1-eu^2)$. Since $u^2 \in {\rm U}_3$, $1-eu^2 \equiv 1-e (\rm{mod}\ 8)$. This implies that $1-eu^2 \notin {\rm N(L^*)}$ and thus $x^2-e \notin \rm{N(L^*)}$. Therefore $x \in \rm{M}$ and $\chi(x)=(1,1)$. \\ 

\noindent {\bf Case $v(e)=2$}:\\
 Suppose $x \in \rm{U}_2 \backslash {\rm U}_{3}$. Then $x \not\in \rm{N(L^*)}$ by Cor.~\ref{important}(iii). Now, $x^2 \in {\rm U}_{3}$ and $v(e) = 2$ implies that $x^2 - e \in {\rm U}_{2} \backslash {\rm U}_{3}$. Therefore $x^2 - e \not\in \rm{N(L^*)}$ by Cor.~\ref{important}(iii). Thus, $\chi(x) = (1,1)$. \\


\noindent {\bf Case $v(e)=3$}:  Choose an element $x=2u$ such that $x \notin {\rm N(L^*)}$. Then $x^2-e = 4(u^2-e/4)$. We have, $u^2-e/4 \equiv (1-e/4)(\rm{mod}\ 8)$. Since $s({\rm L}/\Q_2) =2 $, by remark \ref{rmk:ex}, we get that $v(d)= 1$. Thus by Lemma \ref{e1d1}(ii), at least one of the elements $-1,e,(1-e/4)$ does not belong to ${\rm N(L^*)}$. By our assumption, $-1,e \in {\rm N(L^*)}$. Thus $1-e/4 \notin {\rm N(L^*)}$ which implies that $x^2-e \notin {\rm N(L^*)}$. Therefore $x \in {\rm M}$ and $\chi(x)=(1,1)$. \\ 

\noindent  This completes the proof of the theorem \ref{p2}.

\end{proof}


\section{The irreducible cubic case}\label{section:irreducible}

Let ${\rm K}$ be any finite extension of $\Q_p$. In this section ${\rm X}$ will denote a smooth projective surface $\rm{K}$-birational to the surface defined by the equation
$$y^2-dz^2 = f(x)$$  where $f(x) = x^3+ax^2+bx+c$ is an irreducible monic cubic polynomial with coefficients in ${\rm K}$. In this section we prove Theorem \ref{thm:irreducible} which says that ${\rm A_0(X)_0}=\{0\}$. As mentioned in section \ref{intro}, $d \in \rm{K^*}^2$ implies that ${\rm X}$ is ${\rm K}$-birational to $\P^2$ and hence ${\rm A_0(X)_0}=\{0\}$. Thus we may assume $d\notin {\rm K}^{*2}$. Let $\ \alpha_1,\  \alpha_2,\  \alpha_3\ $ be the roots of $f(x)$ in an algebraic closure of ${\rm K}$. Let ${\rm E}_i = {\rm K}(\alpha_i)$ and ${\rm L = K}(\sqrt{d})$.

\subsection*{The method of computation of ${\rm A_0(X)_0}$}\label{sec3}
 Let ${\rm M}$ = $\{ x \in {\rm K} | x^3+ax^2+bx+c \in \rm{N_{L/K}L^*}\}$. As in  section \ref{section:computation}, results of \cite{Col} reduce the problem of computing the Chow group to the determination of the image of $\rm{M}$ under the map,
$$\chi: \rm{M} \rightarrow \rm{E_1^*/N_{LE_1/E_1}({LE_1}^*}) \times \rm{E_2^*/ N_{LE_2/E_2}(LE_2^*}) \times  \rm{E_3^*/ N_{LE_3/E_3}(LE_3^*})$$
$$x \mapsto (x-\alpha_1, x-\alpha_2, x-\alpha_3) $$

\begin{lemma}\label{lem5} $x-\alpha_i \in {\rm N_{LE_{\it i}/E_{\it i}}({LE_{\it i}}^*})$ if and only if $\ x^3+ax^2+bx+c \in \rm{N_{L/K}L^*}$.
\end{lemma}
\begin{proof}As in the proof of Lemma \ref{lem2}, class field theory implies that $\ x-\alpha_i \in {\rm N_{LE_{\it i}/E_{\it i}}({LE_{\it i}}^*)}$ if and only if ${\rm N_{E_{\it i}/K}({\it x}-\alpha_{\it i}) \in N_{L/K}L^*}$. Since ${\rm N_{E_{\it i}/K}}({\it x}-\alpha_{\it i}) = x^3+ax^2+bx+c$ the lemma follows.

\end{proof}
 
\begin{proof}[Proof of \ref{thm:irreducible}] By definition of ${\rm M}$,  $ x \in {\rm M}$ implies $x^3+ax^2+bx+c \in \rm{N_{L/K}L^*}$. Therefore by Lemma \ref{lem5} we get $x-\alpha_i \in {\rm N_{LE_{\it i}/E_{\it i}}({LE_{\it i}}^*}) \ \text{for} \ i=1,2,3$. Thus $\chi(x) = (0,0,0)$. Thus ${\rm A_0(X)_0}\cong \{0\}$.

\end{proof}


\section{ The Global Case}\label{sec:global} 

\noindent Let ${\rm K}$ be a number field. Let ${\rm X}$ be any smooth projective surface ${\rm K}$-birational to the surface given by the equation $$ y^2-dz^2=x(x^2-e)$$ where $d \notin {\rm K}^{*2}$. Let ${\rm K}_v$ be the completion of ${\rm K}$ at $v$ and ${\rm X}_v= {\rm X}\times_{\rm K}{\rm K}_v$. 

\vspace{.2cm}
\noindent By the result of Bloch (\cite{bloch}, Theorem(0.4)), ${\rm A_0(X_{\it v})_0} = 0$ for almost all places $v$ of ${\rm K}$. We can also observe this directly as follows. At almost all places of ${\rm K}$, ${\rm K}_v(\sqrt d)$ and ${\rm K}_v(\sqrt e)$ are unramified extensions of ${\rm K}_v$. If $d \in {{\rm K}_v^*}^2$, then we have already seen in the section \ref{intro} that ${\rm A_0(X_{\it v})_0} =0$. If $e \in {{\rm K}_v^*}^2$ and  $v$ is not a place lying above the prime ideal $(2)$ then, ${\sqrt e},-{\sqrt e},2{\sqrt e}$ are units. Thus by {\rm (\cite{coray}, Prop.~4.7)}, ${\rm A_0(X_{\it v})_0} = 0$. If neither $e$, nor $d$ is in ${{\rm K}_v^*}^2$, then ${\rm K}_v(\sqrt d)$ and ${\rm K}_v(\sqrt e)$ are both  unramified quadratic extensions and thus isomorphic extensions of ${\rm K}_v$. In this case, by Theorem \ref{L=E} of this paper we get that ${\rm A_0(X_{\it v})_0} = 0$. \\

\noindent Let us recall briefly how the results of Colliot-Th\'el\`ene, Sansuc and Swinnerton-Dyer {\rm(\cite{sansuc-dyer})},  Colliot-Th\'el\`ene  and Sansuc \cite{Col} (see Salberger \cite{salberger} for more general statement valid for all conic bundles over $\P^1_{\rm K}$ ),   which  allow one to compute ${\rm A_0(X)_0}$ from the knowledge of ${\rm A}_0({\rm X}_v)_0$. For each place $v$ of ${\rm K}$, we have a map $ {\rm A_0(X)_0} \to {\rm A_0(X_{\it v})_0}$  and hence, a diagonal map 
$$ \delta : {\rm A_0(X)_0} \to \prod_v {\rm A_0(X_{\it v})_0}.$$ 
\noindent As we have seen above, ${\rm A_0(X_{\it v})_0} = 0$ for almost all places $v$, the target of $\delta$ is the same as $\bigoplus_v{\rm A_0(X_{\it v})_0}$. The exactness of the following sequence is proved in {\rm (\cite{sansuc-dyer1}, Sec.~8)}
$$ 0 \to \Sha^1({\rm K}, {\rm S}) \to {\rm A_0(X)_0} \stackrel{\delta}\to \bigoplus_v {\rm A_0(X_{\it v})_0} \to {\rm Hom(H^1(K,\hat{S}), \Q/\Z)}$$
\noindent  in which $\hat{S}$ = ${\rm Pic (\overline{X})}$, a free $\Z$-module of finite rank with a ${\rm Gal(\overline{K}/K)}$-action, and ${\rm S}$ is the ${\rm K}$-torus dual to $\hat{S}$. The exactness of this sequence reduces the computation of ${\rm A_0(X)_0}$ to the local problem of computing ${\rm A_0(X_{\it v})_0}$. 

\vspace{.2cm}
\noindent Let us indicate how the results of this paper contribute to the solution of this local problem, atleast when ${\rm K} = \Q$. Let ${\rm X}$ be a ${\rm Ch\hat{a}telet}$ surface as above. When $v$ is a finite place of $\Q$, ${\rm A}_0({\rm X}_v)_0$ can be calculated using results of this paper when $x^2-e$ remains irreducible and \cite{dalawat1},\cite{dalawat2} when $x^2-e$ splits into  linear factors over $v$. One now needs to do these calculations when $v$ is the real place of $\Q$, i.e. ${\rm K}_v = \R$. If $x^2-e$ remains irreducible over $\R$, then either ${\rm L}/\R$ is a trivial extension or the extensions ${\rm L}$ and ${\rm E}$ become isomorphic. Thus in both these cases, ${\rm A_0}({\rm X}_v)_0=\{0\}$. When $x^2-e$ is reducible, we use results of \cite{ischebeck} to calculate ${\rm A_0(X_{\it v})_0}$.


\vspace{.2cm}
\noindent Supriya Pisolkar\\
\noindent Harish-Chandra Research Institute\\
\noindent Chhatnag Road, Jhunsi\\
\noindent Allahabad -211019, India.\\
\noindent {\it Email - supriya@mri.ernet.in}

\end{document}